\documentclass[12pt]{amsart}

\usepackage[T2A]{fontenc}
\usepackage[utf8]{inputenc}
\usepackage[english]{babel}

\usepackage{amsmath}
\usepackage{amssymb}
\usepackage{amscd}
\usepackage{amsfonts,amsmath,amsthm}

\usepackage[matrix,arrow,curve]{xy}
\usepackage{xfrac}
\usepackage{color}

\newtheorem{thm}{Theorem}
\newtheorem{Ex}{Example}

\newtheorem{prop}{Proposition}

\newtheorem{lemma}{Lemma}
\newtheorem*{main lemma}{The main lemma}

\newtheorem{definition}{Definition}

\title[Diagonal complexes:  holes  and  involutions]{Diagonal complexes for  {surfaces of finite type} and surfaces with involution }

\author{Joseph Gordon, Gaiane Panina}

\address{  J. Gordon:Chebyshev Laboratory, St. Petersburg State University, St. Petersburg, Russia
International Laboratory of Game Theory and Decision Making, National Research University Higher School of Economics, St. Petersburg, Russia  joseph-gordon@yandex.ru; \ \ G. Panina: St. Petersburg department of Steklov institute of mathematics, Mathematics and Mechanics Faculty, St. Petersburg State University
 gaiane-panina@rambler.ru}

\begin{document}

\begin{abstract} Following the ideas of \cite{GorPan}, we study two related constructions:

---  The  \textit{ diagonal complex}  $\mathcal{D}$ and its barycentric subdivision $\mathcal{BD}$  related to a   oriented surface of a finite type  $F$ equipped with a number of  labeled marked points. This time, unlike \cite{GorPan}, we allow boundary components without marked points, called \textit{holes}.

---  The \textit{symmetric diagonal complex} $\mathcal{D}^{inv}$ and its barycentric subdivision $\mathcal{BD}^{inv}$  related to a  \textit{symmetric}(=with an involution) oriented surface $F$ equipped with a number of (symmetrically placed) labeled marked points.

The symmetric complex  is shown to be homotopy equivalent to the complex of a  surface  obtained by  ''taking a half''  of the initial symmetric surface.

\end{abstract}

\maketitle

\section{Introduction}\label{SectIntro}
Let us start with two elementary motivations:

(*) It is known that the poset of the collections of non-crossing diagonals in an $n$-gon is combinatorially isomorphic to a convex polytope called \textit{associahedron} (see Appendix B).

(**) The poset of the collections of non-crossing diagonals in an  $n$-gon with a hole is combinatorially isomorphic to a convex polytope called \textit{cyclohedron} (see Example \ref{ExCyclohedron}).

A natural question which proved to be meaningful is: \textit{what happens if one replaces the $n$-gon by an arbitrary closed surface equipped with a number of labeled marked points (=vertices)?}

\medskip

Two more sophisticated motivations are:

(***)  Curve complexes (or arc complexes) exist in the literature in different frameworks and settings, see \cite{Thur} and \cite{Harvey} for pioneer papers.    Oversimplifying, the basic idea is to take a (possibly bordered) surface with a finite set of labeled distinguished points, and to associate a complex with the ground set (that is, the set of vertices) equal to  homotopy classes of either closed curves, or curves with endpoints in the distinguished set, or both, as in  \cite{KorPap}. Simplices correspond to non-intersecting representatives of the homotopy classes.
The  mapping class group has a natural subgroup acting  on the complex, so it makes sense to take the quotient space.

(****) Combinatorial models in the moduli spaces theory is a classical object which is naturally linked to the present research.

\bigskip

The motivation (*) has led  to
  \cite{GorPan}, where we introduced and studied the complex of pairwise non-intersecting diagonals on an oriented surface equipped with $n$ marked points.
If the surface of genus $g$ is closed, the complex is homotopy equivalent to the space of metric ribbon graphs $RG_{g,n}^{met}$, or, equivalently, to the decorated moduli space $\widetilde{\mathcal{M}}_{g,n}= \mathcal{M}_{g,n}\times\mathbb{R}_+^n$. For bordered surfaces, we proved in \cite{GorPan} the following: (1) Contraction of a boundary edge does not change the homotopy type of the support of the complex. (2)
Contraction of a boundary component to a new marked point yields a forgetful map between two diagonal complexes which is homotopy equivalent to the Kontsevich's tautological circle bundle $L_i$. (3) In the same way, contraction of several boundary components corresponds to Whitney sum of the tautological bundles.

It is important that the paper \cite{GorPan} deals with the case when each of the boundary components contains at least one marked point.
Motivation (**) suggests us to relax this condition.
In the present paper we allow boundary components without marked points, or \textit{holes}. Filling in a hole
 gives rise to  a bundle whose fiber is a surface $\overline{F}$ obtained from $F$ by eliminating a small disk around each of the marked points, see Section \ref{SecPuncMain}.

In Section \ref{SectSymm} we consider
 symmetric surfaces (that is, surfaces with a distinguished involution), symmetric diagonal arrangements, and associated diagonal complex. The latter is shown to be homotopy equivalent to the diagonal complex of   surface  obtained by ''halving'' of the initial surface.
Oversimplifying, the relation reads as follows:  cut the surface $F$ through the "symmetry axis" and take one half of $F$. The cut line turns to a hole.   This naive surgery  leads to a map between diagonal complexes  which is a homotopy equivalence, see Theorem \ref{ThmHalf}.

\subsection*{Acknowledgement} This research is
 supported by the Russian Science Foundation under grant 16-11-10039.

 We are also indebted to Peter Zograf and Max Karev for useful remarks.

\section{Diagonal complex: construction and introductory examples}\label{SectMainConstr}

Assume that  an oriented surface $F$ of genus $g$ with $b+f$ labeled boundary components $B_i$ is fixed.  We mark  $n$ distinct labeled points on $F$ not lying on the boundary. Besides, for each $i=1,..,b$  we fix $n_i>0$ distinct labeled points on  the  boundary component $B_i$.

So, there are $f$ boundary components without marked points. We call them \emph{holes}.

Let us stress  once again, that holes were not allowed in \cite{GorPan}.

  We assume that $F$ can be triangulated (decomposed into triangles, possibly with holes) with vertices at the marked points.

 Altogether we have $N=n+\sum_{i=1}^b n_i$ marked points; let us call them  \textit{vertices} of $F$. The vertices not lying on the boundary are called \textit{free vertices}. The vertices that lie on the boundary split the boundary components into \textit{edges}.

A \textit{pure diffeomorphism}  $F \rightarrow F$ is an orientation preserving  diffeomorphism which maps marked points
to marked points,  and preserves the labeling of marked points and holes.  Therefore, a pure diffeomorphism
maps each boundary component to itself. The \textit{pure mapping class group} $PMC(F)$ is the group of isotopy classes of pure diffeomorphisms.

A \textit{diagonal}  is a simple (that is, not self-intersecting) smooth curve  $d$  on $F$ whose endpoints are some of the (possibly the same) vertices  such that\begin{enumerate}
 \item $d$ contains no vertices (except for the endpoints).
                     \item $d$ does not intersect the boundary (except for its endpoints).
                      \item $d$ is not  homotopic to an edge of the boundary.

                      Here and in the sequel, we mean homotopy with fixed endpoints in the complement  of the vertices $F \setminus \  Vert$.  In other words, a homotopy never hits a vertex.
                     \item $d$ is non-contractible.

                   \end{enumerate}

An \textit{admissible diagonal arrangement} (or an \textit{admissible arrangement}, for short) is a collection of diagonals $\{d_j\}$ with the properties:

\begin{enumerate}
\item Each free vertex is an endpoint of some diagonal.
  \item No two diagonals intersect (except for their endpoints).
  \item No two diagonals are homotopic.
  \item The complement of the arrangement and the boundary components $(F \setminus \bigcup d_j) \setminus \bigcup B_i$  is a disjoint union of open (possibly with holes) disks.
We allow any number of holes in a single disc.
\end{enumerate}

 We say that a tuple $(g,b,n,f)$ is \textit{stable
} if no admissible arrangement has  a non-trivial automorphism (that is, each pure diffeomorphism which maps an arrangement to itself, maps each germ of each of the diagonals $d_i$ to itself).

 Tuples with $b>1$ are stable since a boundary component allows to set a linear ordering on the germs of diagonals emanating from each of its vertices.
It is known\footnote{This follows from Lefschetz fixed point theorem, as explained  by Bruno Joyal in personal communications.} that any tuple with $n>2g+2$ is stable.
Throughout the paper we assume that \textbf{all  the tuples are stable}.

\begin{definition}\label{DefEquivArr}
Two arrangements  $A_1$ and $A_2$ are \textit{strongly equivalent}  whenever there exists a homotopy taking $A_1$ to $A_2$.

Two arrangements  $A_1$ and $A_2$ are \textit{weakly equivalent}  whenever there exists a composition of a homotopy and a pure diffeomorphism of $F$ which maps bijectively $A_1$ to $A_2$.
\end{definition}

\subsection*{Poset $\widetilde{D}$ and cell complex $\widetilde{\mathcal{D}}$.}

Strong equivalence classes of admissible arrangements are partially ordered by reversed inclusion: we say that $A_1 \leq A_2$ if there exists a homotopy that takes the arrangement $A_2$ to some subarrangement of $A_1$.

Thus  for the data  $ (g, b , n,f; n_1,...,n_b)$ we have the posets of all strong equivalence classes of admissible arrangements  $\widetilde{D}={\widetilde{D}}_{g, b , n,f; n_1, \ldots,n_b}$.

\begin{thm} \label{ThmOldPunc} (\cite{Panina},\cite{Pen1}). If $F$ is a polygon with $n_1>0$  marked points on its boundary, $f$ holes, and no free marked points, then  the poset $\widetilde{D}={\widetilde{D}}_{0, 1 , 0,f; n_1}$ is a combinatorial ball $B^{n_1+2f-3}$. The cellulation has the unique biggest cell that corresponds to the empty diagonal arrangement.\qed
\end{thm}

The following example (and its generalizations) is well-known  and is used in the cluster algebras world \cite{FomZel}, \cite{FraScSo}. However we put it here with a proof for the sake of completeness.

\begin{figure}[h]
\centering \includegraphics[width=6 cm]{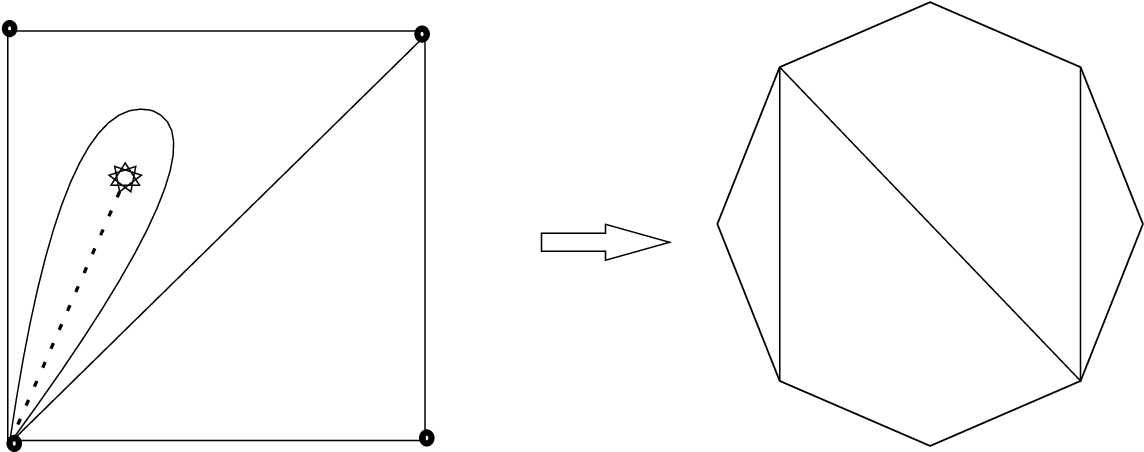}
\caption{Building a symmetric diagonal arrangement. The dashed line denotes the cut.}\label{cyclohedr}
\end{figure}

\begin{Ex} \label{ExCyclohedron} The complex  $\mathcal{D}_{0,1,0,1;n}$ (that is, the diagonal complex associated with an $n$-gon with a hole) is combinatorially isomorphic to the cyclohedron.
\end{Ex}
Proof. Given an admissible arrangement, cut the polygon by a path   connecting the hole with the boundary of the polygon. We assume that the cut does not cross the diagonals. Take the copy of the polygon with the same arrangement and with the same cut, and glue the two copies together. We get a    $2n$-gon together with a centrally symmetric diagonal arrangement, see Fig. \ref{cyclohedr}.  This construction can be reversed, and therefore establishes a combinatorial isomorphism with the cyclohedron (see Appendix B). \qed

\medskip

Theorem \ref{ThmOldPunc} implies that
the poset $\widetilde{D}$ can be realized as the poset of some (uniquely defined)  regular\footnote{
A cell complex $K$ is \textit{regular} if each  $k$-dimensional cell $c$ is attached to some  subcomplex of
the $(k-1)$-skeleton of $K$ via a bijective mapping on $\partial c$.} cell complex $\widetilde{\mathcal{D}}$.
Indeed, let us build up  $\widetilde{\mathcal{D}}$  starting from the cells of maximal dimension. Each such cell corresponds to cutting of the surface $F$ into a single polygon with $f$ holes. Adding more diagonals reduces to Theorem \ref{ThmOldPunc}. In other words,   $\widetilde{\mathcal{D}}$ is a patch of combinatorial balls that arise in the theorem.

 For the most examples, $\widetilde{\mathcal{D}}$ has infinitely many cells. Our goal is to factorize $\widetilde{\mathcal{D}}$ by the action of the pure mapping class group. For this purpose consider the defined below barycentric subdivision of $\widetilde{\mathcal{D}}$.

\subsection*{Poset $\widetilde{BD}$ and cell complex $\widetilde{\mathcal{BD}}$.} We apply now the construction of the \textit{order complex} \cite{Wachs} of a poset, which yields the barycentric subdivision.
  Each  element of the  poset ${\widetilde{BD}}_{g, b , n,f; n_1,... ,n_b}$ is (the strong equivalence class of) some admissible arrangement $A=\{d_1,...,d_m\} $ with a linearly ordered partition $A=\bigsqcup S_i$   into some non-empty sets  $S_i$  such  that the first set $S_1$ in the partition  is an admissible arrangement.


The partial order on $\widetilde{BD}$ is generated by the following rule:
 \newline $(S_1,...,S_p)\leq (S'_1,...,S'_{p'})$ whenever one of the two conditions holds:
\begin{enumerate}
  \item We have one and the same arrangement  $A$, and $(S'_1,...,S'_{p'})$ is an order preserving refinement of $(S_1,...,S_p)$.
  \item $p\leq p'$, and for all $i=1,2,...,p$, we have  $S_i=S'_i$. That is, $(S_1,...,S_p)$ is obtained from $(S'_1,...,S'_{p'})$ by removal
  $S'_{p+1},...,S'_{p'}$.
\end{enumerate}

Let us look at the incidence rules in more details. Given $(S_1,...,S_p)$, to list all the elements of $\widetilde{BD}$ that are smaller than $(S_1,...,S_p)$  one has
(1) to eliminate some (but not all!) of $S_i$ from the end of the string, and (2) to replace some consecutive collections of sets by their unions.

\medskip

Examples:$$(\{d_5,d_2\},\{d_3\},\{d_1,d_6\},\{d_4\},\{d_7\},\{d_8\})\ \  >\ \ (\{d_5,d_2\},\{d_3,d_1,d_6\},\{d_4,d_7\}).$$
$$(\{d_5,d_2\},\{d_3\},\{d_1,d_6\},\{d_4\},\{d_7\},\{d_8\})\ \  >\ \ (\{d_5,d_2\},\{d_3\},\{d_1,d_6\},\{d_4\},\{d_7\}).$$
$$(\{d_5,d_2\},\{d_3\},\{d_1,d_6\},\{d_4\},\{d_7\},\{d_8\})\ \  >\ \ (\{d_5,d_2\},\{d_3\},\{d_1,d_6\},\{d_4\},\{d_7,d_8\}).$$

Minimal elements of $\widetilde{BD}$ correspond to admissible arrangements. Maximal elements correspond to maximal arrangments  $A$ together with some minimal admissible subarrangement $A'\subset A$ and a linear ordering on the set $A\setminus A'$.

{By construction}, the complex  $\widetilde{\mathcal{BD}}$  is combinatorialy isomorphic to the barycentric subdivision of $\widetilde{\mathcal{\mathcal{D}}}$.

We are mainly interested in the  quotient complex:

\begin{definition}
  For a fixed data $(g, b , n,f; n_1,... ,n_b)$, the diagonal  complex $\mathcal{BD}_{g, b , n,f; n_1,... ,n_b}$ is defined as

  $$\mathcal{BD}=\mathcal{BD}_{g, b , n,f; n_1,... ,n_b}:=\mathcal{\widetilde{BD}}_{g, b , n,f; n_1,... ,n_b}/PMC(F).$$
  We define also
   $$\mathcal{D}=\mathcal{{D}}_{g, b , n,f; n_1,... ,n_b}:=\mathcal{\widetilde{D}}_{g, b , n,f; n_1,... ,n_b}/PMC(F).$$
\end{definition}
 Alternative definition reads as:

 \begin{definition}
   Each cell of the complex ${{\mathcal{BD}}}_{g, b , n,f; n_1,... ,n_b}$ is labeled by the weak equivalence class of some admissible arrangement $A=\{d_1,...,d_m\} $ with a linearly ordered partition $A=\bigsqcup S_i$   into some non-empty sets  $S_i$  such  that the first set $S_1$   is an admissible arrangement.

The incidence rules  are the same as the above rules for the complex $\widetilde{\mathcal{BD}}$.
 \end{definition}

\begin{prop} The cell complex $\mathcal{BD}$ is regular. Its cells are combinatorial simplices.

\end{prop}
Proof.   If $(S_1,...,S_r)\leq (S'_1,...,S'_{r'})$ then there exists a unique (up to isotopy)  order-preserving  pure diffeomorphism of $F$ which  embeds \newline  $A= S_1\cup ...\cup S_r$ in  $A'=S'_1\cup ...\cup S'_{r'}$.
Indeed,   If $S_1=S_1'$,  the arrangement $S_1$  maps identically to itself since it has no automorphisms by stability assumption. The rest of the diagonals are diagonals in polygons, and are uniquely defined by their endpoints. Assume that $S_1 \subset S_1'$.
For the rest of the cases it suffices to take $A=S_1$, $A'=A =S_1'\bigsqcup S_2'$. If $A$ embeds in $A'$ in different ways,
then $A$ has a non-trivial isomorphism, which contradicts stability assumption.
\qed

\medskip

\textbf{Remark.}
A reader may imagine
 each (combinatorial) simplex in ${\mathcal{BD}}$ as a (Euclidean) equilateral simplex and to
define the \textit{support}, or \textit{geometric realization} of the complex $| {\mathcal{BD}}|= |{\mathcal{D}}|$  as the patch of these simplices.

\section{Diagonal complexes related to holed surfaces: main theorems}\label{SecPuncMain}

Filling in a hole gives rise to
a natural forgetful projection$$\pi: \mathcal{BD}_{g, b , n,f+1; n_1,... ,n_b} \rightarrow  \mathcal{BD}_{g, b , n,f; n_1,... ,n_b}.$$

It is defined as follows.  An element of $\mathcal{BD}_{g, b , n,f+1; n_1,... ,n_b}$  corresponds to some   admissible arrangement together with a partition
\((S_1,...,S_r)\). Fill in the hole $f+1$. We obtain a  collection of diagonals on the surface with $f$ holes. Some of the diagonals may become either contractible or homotopic to an edge of \(F\). Eliminate them.
Some of the diagonals may become pairwise homotopy equivalent. In each class we leave exactly one that belongs to \(S_i\) with the smallest index \(i\).
Eventually some of the sets \(S_i\) may become empty in the process. Eliminate all the empty sets keeping the order of the rest. We obtain an element from $BD_{g, b , n,f; n_1,... ,n_b}$. It is easy to check that $A<A'$ implies $\pi(A)\leq \pi(A')$, so the map is indeed a poset morphism.

 The poset morphism extends to a   piecewise linear map (we denote it by the same letter $\pi$).

We shall need the new surface $\overline{F}$ which is obtained from  $F$  by filling in all the holes,  and replacing each of the free marked points $v_i$ by a hole $H_i$.

\begin{thm}\label{ThmContrFree} (1) The above defined forgetful projection

$$\pi: \mathcal{BD}_{g, b , n,f+1; n_1,... ,n_b} \rightarrow  \mathcal{BD}_{g, b , n,f; n_1,... ,n_b}$$
is  homotopy equivalent to a locally trivial bundle over $\mathcal{BD}_{g, b , n,f; n_1,... ,n_b}$  whose fibers are homeomorphic to the (above defined) surface $\overline{F}$.

(2) Each of the holes $H_i$ gives rise to a circle bundle over $\mathcal{BD}_{g, b , n,f; n_1,... ,n_b}$. It is  isomorphic to the tautological
circle bundle $L_i$.
\end{thm}
Proof.
Let us examine the preimage  $\pi^{-1}(x)$  of an inner point $x$  of a simplex of $\mathcal{BD}_{g, b , n,f; n_1,... ,n_b}$ labeled by $(S_1,...,S_r)$. The label corresponds to an arrangement $A=\bigcup S_i$ on $F$ with $f$ holes. We shall show that $\pi^{-1}(x)$
 is a cell complex homeomorphic to $\overline{F}$.

The explicit construction  consists of two steps.  On the  first step we understand what diagonals  are added to $A$  in the labels of the cells
that intersect $\pi^{-1}(x)$.

On the second step we analyze the partition on the extended set $A$.

\medskip

\textbf{Step 1. }Analyze the ''new'' curves in the preimage.

Remove all the free boundary components (or holes) from $F$. Now
the arrangement $A$ cuts  $F$ into some  polygons.  A \textit{corner} is a vertex with two germs of incident edges $g_1$ and $g_2$ such that there are no other germs between  $g_1$ and $g_2$. For each of the corners, we blow up its vertex, that is, replace it by an extra edge, as is shown in Fig. \ref{Corners}. So each free vertex turns to a new boundary component. We get a two-dimensional cell complex $\overline{\mathcal{F}}(A)$ homeomorphic to $\overline{F}$.

Each cell $\sigma$ of $\overline{\mathcal{F}}(A)$  gives rise to a  new diagonal arrangement $A(x, \sigma)\supset A$  and specifies the place where the boundary component $P_{f+1}$ should be inserted.  This is illustrated in Figures \ref{FigPreimage}  and \ref{Fig3}, and is described by the following rules:

\begin{figure}[h]
\centering \includegraphics[width=13 cm]{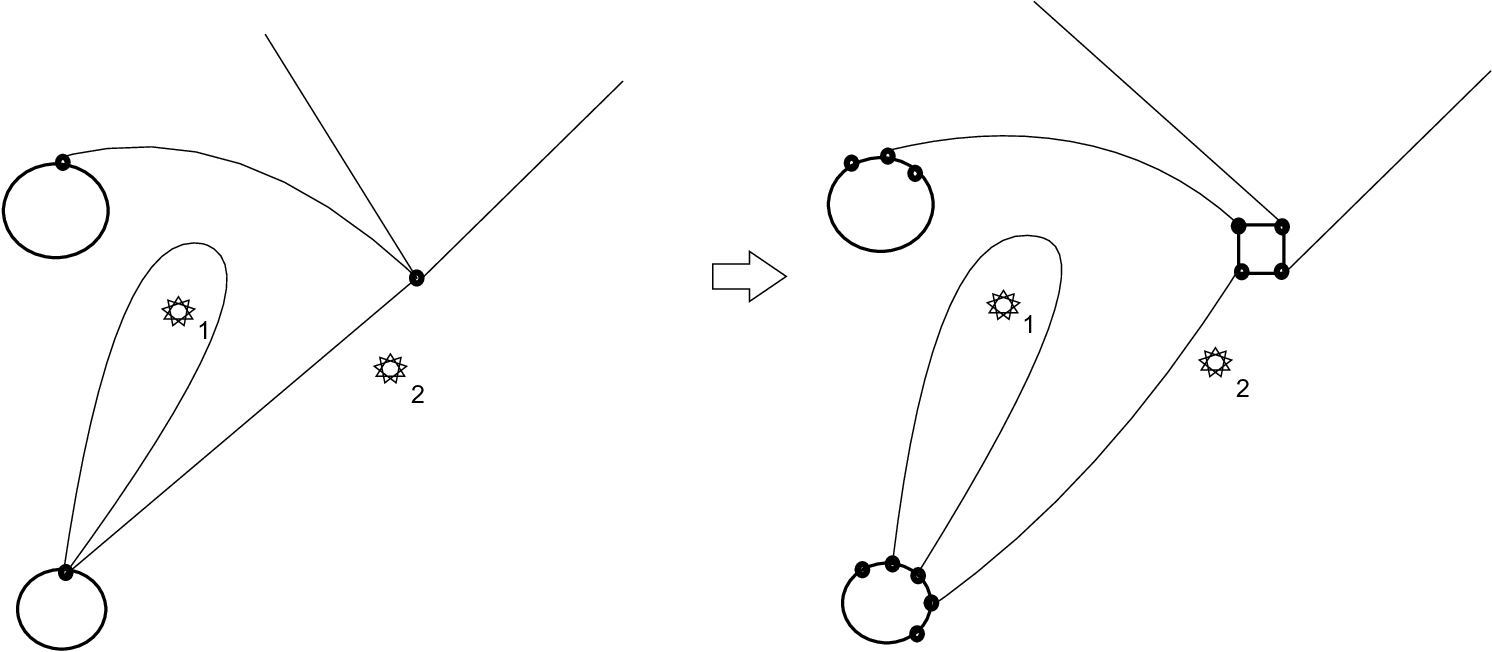}
\caption{Blowing up all the corners. Here we depict a fragment of a bordered surface $F$. Bold circles are the boundary components with marked points, small stars denote the holes.}\label{Corners}
\end{figure}

\begin{figure}[h]
\centering \includegraphics[width=13 cm]{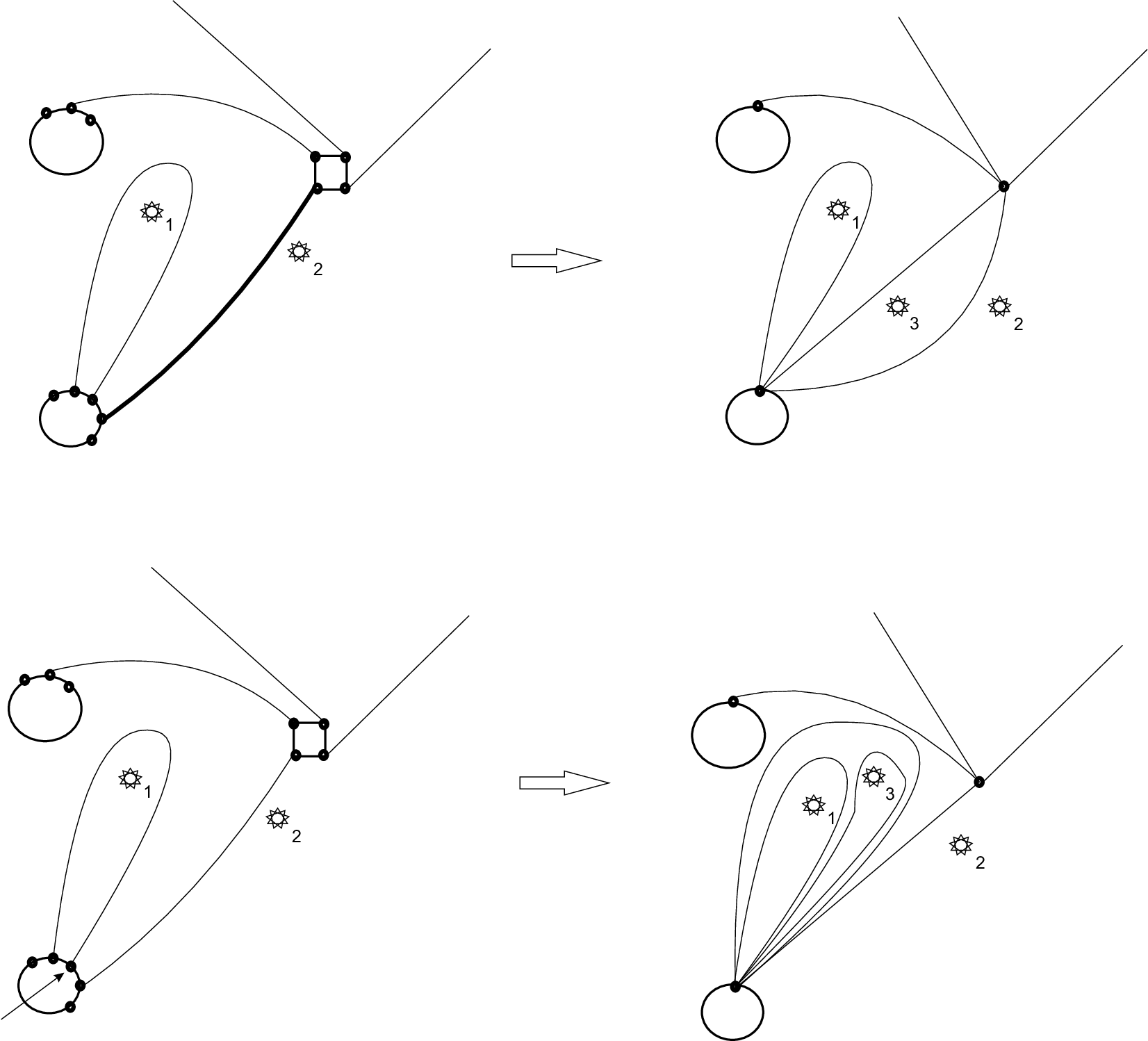}
\caption{Here we depict arrangements that correspond to the bold edge   and to the vertex (indicated by the small arrow) of $\mathcal{F}$.}\label{FigPreimage}
\end{figure}

\begin{figure}[h]
\centering \includegraphics[width=9 cm]{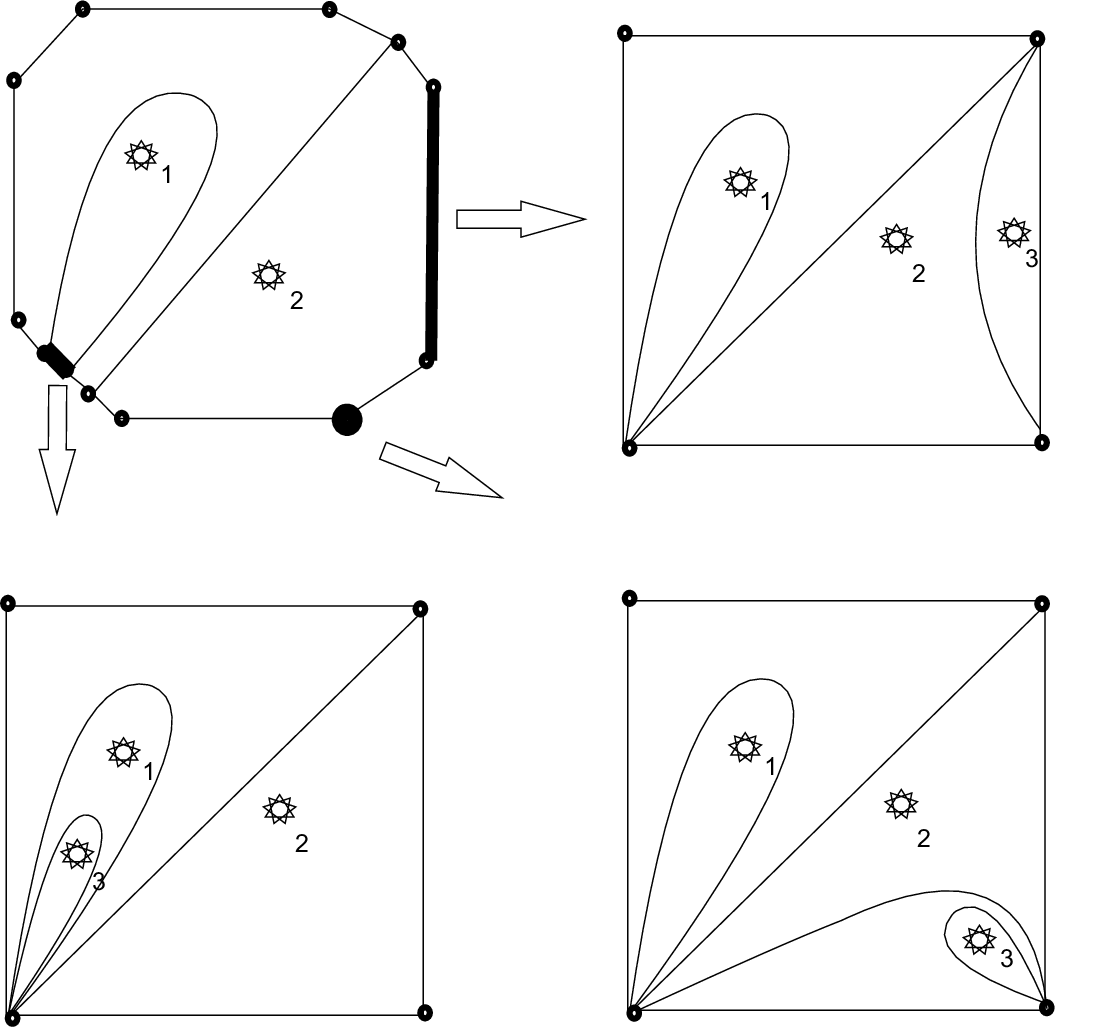}
\caption{One more illustration. Here we depict arrangements that correspond to a vertex and two edges.}\label{Fig3}
\end{figure}

\begin{enumerate}
  \item If $\sigma$ is a $2$-dimensional cell,  keep the arrangement  $A$ as it is, and add the new hole $P_{f+1}$ in the cell $\sigma$.
  \item If the cell $\sigma$ is an edge $e$ of $F$,   add one more diagonal which is parallel to $e$, and add  the new hole $P_{f+1}$  between the two copies of $e$.
  \item  If the cell $\sigma$  is a blow-up of one of the corners, add $P_{f+1}$ in a cell which is adjacent to the corner, and add one loop diagonal embracing  $P_{f+1}$. The loop  starts and ends at the corner.
  \item  If the cell is one of the diagonals  $d\in A$,  duplicate $d$, and put  $P_{f+1}$ between $d$ and its copy.

     \item  If the cell is  one of the new vertices, that is, corresponds to a corner and a diagonal (or to a corner and an edge) of $F$,  we combine either    (3) and (2),  or (3) and (4). That is, add both a loop and a double of a diagonal (or a double of an edge).

\end{enumerate}

\medskip

\textbf{Step 2. }Analyze the partition in the preimage.




First, construct a new cell complex $\widehat{\mathcal{F}}=\widehat{\mathcal{F}}(S_1,...,S_r)$  whose support is homeomorphic to $\overline{\mathcal{F}}$.  This complex is more sensitive: it ''knows'' the partition $A=\bigcup S_i$, whereas $\overline{\mathcal{F}}$ depends on $A$ only.
\begin{enumerate}
  \item Start with the complex $\overline{\mathcal{F}}$.
   Replace each  boundary edge  of $\overline{\mathcal{F}}$ by $r+1$ parallel lines. We call the area between the lines the \textit{grid related to the edge}.
  \item Replace each of the internal  edges of $\overline{\mathcal{F}}$, that is, each of diagonals $d\in A$ by $2(r-k)+3$ parallel lines, where $k$ is defined by  condition  $d\in S_k$. The area between the lines is called the \textit{grid related to the edge}.
  \item Each of the vertices is thus replaced by a rectangular grid $r\times (2(r-k)+3)$. Add diagonals to some of the squares as is shown in Figure  \ref{FigGrid}.

      \item  Contract the  strips corresponding to edges and diagonals along their lengths, so each of the strips becomes a segment.  The rectangular grids survive unchanged.

           \item  The complement of the union of all the grid areas is a number of ''\textit{big disks}''.  They bijectively correspond to two-cells of $\overline{\mathcal{F}}$. Contract each of them to a point.

\end{enumerate}

\medskip

\begin{figure}[h]
\centering \includegraphics[width=14cm]{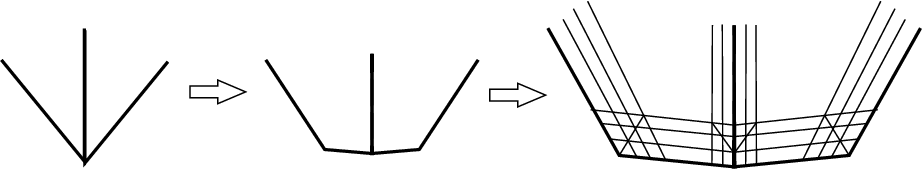}
\caption{These are corresponding fragments of $F$, $\overline{\mathcal{F}}$, and   $\widehat{\mathcal{F}}$. We depict    grids but  do not depict contractions. }\label{FigFGrid}
\end{figure}

\textbf{Remark.} In Figures \ref{FigFGrid},\ref{FigGrid}  we depict a piece of $\widehat{\mathcal{F}}$  \textbf{without }contractions  for the sake of better visualization. Keeping the contractions in mind, one should remember that  the contracted cells have a different dimension.

\medskip

The intersection of $\pi^{-1}(x)$   with the cells of the complex \newline
$\mathcal{BD}= \mathcal{BD}_{g, b , n,f+1; n_1,... ,n_b}$  yields  a structure of a cell complex
on $\pi^{-1}(x)$.
Let us show that the  cell structure is isomorphic to $\widehat{\mathcal{F}}$ by presenting a bijection between the cells of $\widehat{\mathcal{F}}$  and the cells of $\mathcal{BD}$   intersecting $\pi^{-1}(x)$.

We add new diagonals  to $(S_1,...,S_r)$ and $P_{f+1}$ in the same way as we did on Step 1. But this time we  specify the partition on the new set of diagonals. That is, we decide in addition where to put the new diagonals (is any).
Possible options are: we either put them in one of $S_i$, or create separate sets coming right after one of $S_i$. The above described grid tells us
the choice of an option.

\begin{enumerate}
  \item If $\sigma$ is a $0$-dimensional cell  coming from a big disk,   keep the arrangement  $A$, its partition
  $ (S_1,...,S_r)$ and add the new hole $P_{f+1}$ in the corresponding two-cells of $\overline{\mathcal{F}}$.

  \item If the cell $\sigma$ is an edge of a grid related to  one of the edges of $F$, we keep
  $ (S_1,...,S_r)$,  duplicate the edge, and add  $P_{f+1}$  between the edge and its double. The new diagonal (that is, the double of the edge)  we put in one of $S_i$ where $i$ is determined by the number of the line in the grid. The lines of the grid are numbered starting from the edge.

   \item If the cell $\sigma$ is a strip of a grid related to  one of the edges of $F$, we add the same diagonal as above, but now we put it in a singleton coming after $S_i$.   The strips of the grid are numbered  also  starting from the edge.

  \item  We apply  the same strategy  when adding a  loop diagonal embracing the new  $P_{f+1}$.

  \item  If the cell corresponds to one of the diagonals  $d\in S_k$, we duplicate $d$, put  $P_{f+1}$ between $d$ and its copy.
  The new and the old diagonals now are undistinguished. At least one of them should be in   $S_k$, the other one should be either in $S_k$ (this corresponds to the central line of the grid) or to the right of $S_k$.  Altogether we have $(S_1,S_2,...,S_r,\{d\}), (S_1,S_2,...,S_r\cup \{d\}), $ \newline $(S_1,S_2,...,\{d\}, S_r,),...,
(S_1,S_2,...,S_k\cup \{d\},...,S_r),...$ \newline $(S_1,S_2,...,S_k\setminus \{d'\}\cup\{d\},...,S_r\cup \{d'\}),
(S_1,S_2,...,S_k\setminus \{d'\}\cup\{d\},...,S_r, \{d'\})$.

     \item  If the cell  comes from the grid at some vertex   we combine previous construction.

\end{enumerate}

\medskip

\textbf{Example.}

Figure \ref{FigGrid}  depicts a piece of $\widehat{\mathcal{F}}$  with grids.
This fragment corresponds to a (vertical) diagonal $d$ plus two  adjacent  corners, right and left. In this particular example   $r=3$, that is, we have $(S_1,S_2,S_3)$. Here we also have   $k=2$, that is,  $d\in S_2$.  The arrows in the figure point to some cells of $\widehat{\mathcal{F}}$.

\begin{figure}[h]
\centering \includegraphics[width=10cm]{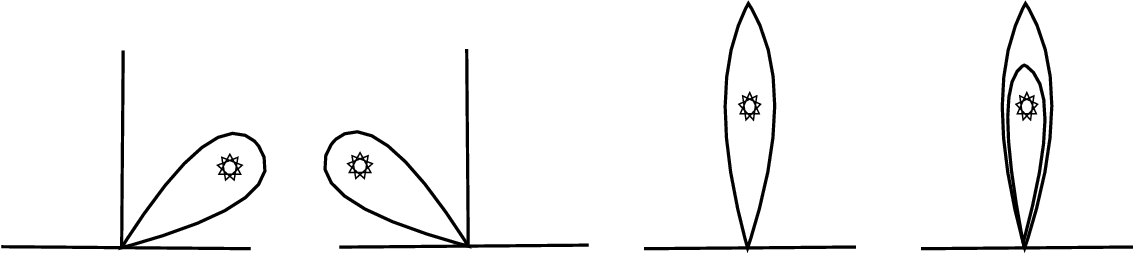}
\caption{Illustration to the Example: right loop  $e_r$,  left loop  $e_l$, duplicated diagonal, loop $+$ duplicated diagonal. On the last figure one sees that right and left loops are indistinguishable whenever the edge is duplicated.  }\label{FigFcomments}
\end{figure}

\begin{figure}[h]
\centering \includegraphics[width=9 cm]{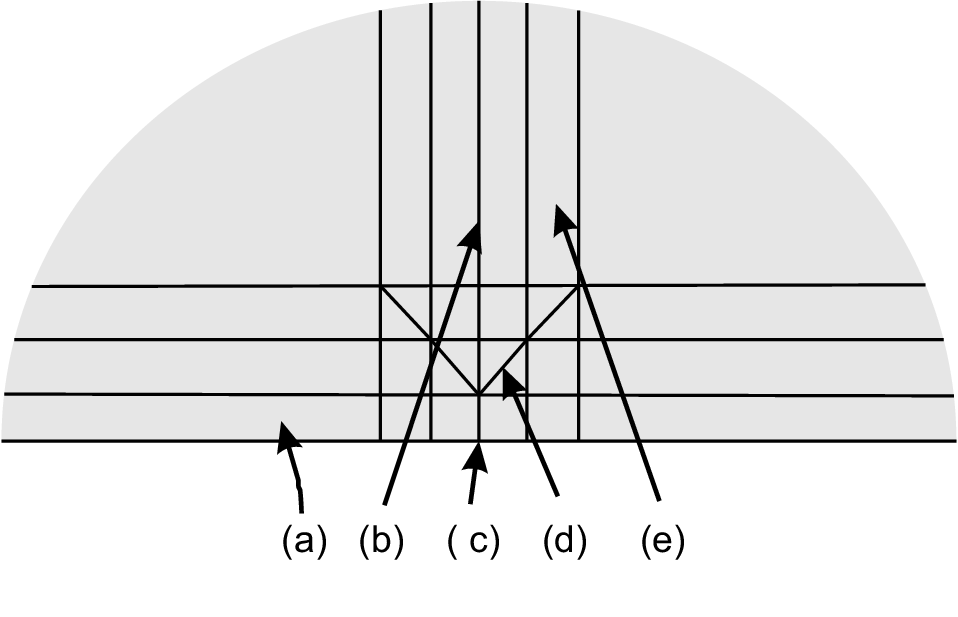}
\caption{This is a fragment  of $\widehat{F}$  with grids but without contractions. The arrows refer to the above example.}\label{FigGrid}
\end{figure}

The cells correspond to adding loops  $e_l$ or  $e_r$, and  duplicating the diagonal $d$, see Figure \ref{FigFcomments}.
The double of the diagonal we denote by $d'$.
If the diagonal is duplicated, the left and right loops are indistinguishable.

We list the labels of the corresponding cells in $\mathcal{BD}$:

(a) corresponds to $(S_1,\{e_l\}, S_2,S_3)$.

(b) corresponds to  $(S_1,S_2\cup \{d'\},S_3)$.

(c)     corresponds to  $(S_1\cup \{e'\},S_2\cup \{ d'\},S_3)$.

(d) corresponds  to $(S_1,S_2, \{e_r\}\cup \{d'\},S_3)$.

(e) corresponds  to $(S_1,S_2,S_3, \{d'\})$.
\qed

\medskip
\textbf{Remark.}

One can rephrase the above theorem as: the triple
 $$\pi: \mathcal{BD}_{g, b , n,f+1; n_1,... ,n_b} \rightarrow  \mathcal{BD}_{g, b , n,f; n_1,... ,n_b}$$

is homotopy equivalent to the universal curve where each of the marked points is replaced by a hole.

\section{Symmetric diagonal complexes: main construction and introductory examples}\label{SectSymm}

We are going to repeat the construction of Section \ref{SectMainConstr}  for symmetric arrangements on a  symmetric surface $F$.
Assume that  an oriented surface $F$ of genus $2g$ with $2b$ labeled boundary components $B_1,...,B_{2b}$ is fixed.

Let us assume that a continuous involution $inv:F\rightarrow F$ is  such that
\textbf{ the set of fixed points is a separating circle} $C^{fix}\subset F$. Such an involution necessarily reverses the orientation.

 We mark  $2n$ distinct labeled points on $F$ not lying on the boundary. Besides, for each $i=1,..,b$  we fix $n_i>0$ distinct labeled points on  each of  the  boundary components $B_i$ {and $B_{b+i}$} assuming   that:

  \begin{enumerate}
    \item no marked point lies on $C^{fix}$;
    \item the involution $inv$ maps marked points to marked points;
		{\item the involution $inv$ maps $B_i$ to $B_{b+i}$ for each $i=1,..,b$.}
  \end{enumerate}

  We assume that $F$ can be tiled by polygons with vertices at the marked points such that: (1) each polygon has at least three vertices, and (2) the tiling is symmetric with respect to $inv$.  For instance, we exclude all  cases like sphere with two marked points.

 Altogether we have $2N=2n+2\sum_{i=1}^b n_i$ marked points; let us call them  \textit{vertices} of $F$.

A \textit{pure diffeomorphism}  $F \rightarrow F$ is an orientation preserving  diffeomorphism which \textbf{commutes with} $inv$ and maps each labeled point to itself.  Therefore, a pure diffeomorphism
maps each boundary component to itself. The \textit{pure mapping class group} $PMC^{inv}(F)$ is the group of isotopy classes of pure diffeomorphisms.

A \textit{diagonal}  is a simple (that is, not self-intersecting) smooth curve  $d$  on $F$ whose endpoints are some of the (possibly the same) vertices  such that\begin{enumerate}
 \item $d$ contains no vertices (except for the endpoints).
                     \item $d$ does not intersect the boundary (except for its endpoints).
                      \item $d$ is not  homotopic to an edge of the boundary.

                     \item $d$ is non-contractible.

                   \end{enumerate}

\begin{lemma}\label{LemmaHomot}
Assume that two homotopic diagonals $d_1$ and $d_2$ connect $v$ and $inv(v)$. Assume also that both $d_1$ and $d_2$ are symmetric: $inv$ maps each of $d_1$ and $d_2$ to itself. Then there exists a symmetric (=commuting with $inv$) homotopy taking $d_1$ to $d_2$.\qed
\end{lemma}

A \textit{symmetric admissible diagonal arrangement} (or a \textit{symmetric admissible arrangement}, for short) is a non-empty collection of diagonals $\{d_j\}$ with the properties:

\begin{enumerate}
\item Each free vertex is an endpoint of some diagonal.
  \item No two diagonals intersect (except for their endpoints).
  \item No two diagonals are homotopic.

  Two remarks are necessary:
                    (a) As in the previous sections,  we mean homotopy with fixed endpoints in the complement  of the vertices $F \setminus \  Vert$.  In other words, a homotopy never hits a vertex.

                    (b) The condition "The homotopy commutes with $inv$" is relaxed due to Lemma \ref{LemmaHomot}.

  \item The complement of the arrangement and the boundary components $(F \setminus \bigcup d_j) \setminus \bigcup B_i$  is a disjoint union of open disks.
      \item $inv$ {takes the arrangement to itself.}
\end{enumerate}

\begin{lemma}\label{Lemma1}
Let $F$ and  $inv$ be as above.
\begin{enumerate}
  \item If a diagonal $d$ belongs to an admissible arrangement, then it intersects $C^{fix}$ at at most one point.
  \item If a disk in $F\setminus \{$diagonals$\}$ is bounded by diagonals  $d_1,...,d_k$, then  all the diagonals taken together intersect $C^{fix}$ at at most two points.
\end{enumerate}
\qed
\end{lemma}

\begin{definition}\label{DefEquivSymArr}
Two symmetric arrangements  $A_1$ and $A_2$ are \textit{strongly equivalent}  whenever there exists a symmetric homotopy taking $A_1$ to $A_2$.

Two symmetric arrangements  $A_1$ and $A_2$ are \textit{weakly equivalent}  whenever there exists a composition of a symmetric  homotopy and a pure diffeomorphism of $F$ which maps bijectively $A_1$ to $A_2$.
\end{definition}

\medskip

We assume that \textbf{all  the tuples are stable}: no admissible arrangement has  a non-trivial  automorphism (that is, each pure diffeomorphism which maps an arrangement to itself, maps each germ of each of $d_i$ to itself).

\medskip

\subsection*{Poset $\widetilde{D}^{inv}$ and cell complex $\widetilde{\mathcal{D}}^{inv}$.}

Strong equivalence classes of symmetric admissible arrangements are partially ordered by reversed inclusion: we say that $A_1 \leq A_2$ if there exists a symmetric homotopy that takes the arrangement $A_2$ to some symmetric subarrangement of $A_1$.

Thus  for the data  $ (g, b , n; n_1,...,n_b)$ we have the poset of all strong equivalence classes of symmetric admissible arrangements  \newline $\widetilde{D}^{inv}={\widetilde{D}}^{inv}_{g, b , n; n_1, \ldots,n_b}$.

Minimal elements of the poset correspond to maximal (by inclusion) symmetric admissible arrangements, that is, to symmetric cuts of $F$
into triangles and quadrilaterals  such that each quadrilateral is symmetric under $inv$.

Maximal elements of the poset correspond to minimal symmetric admissible arrangements, that is,  to symmetric cuts of $F$
into a single disc.

\begin{prop}\label{ExSymmAssoc} Take a planar regular $2k$-gon with a fixed symmetry axis which contains no vertices. Consider the poset of collections of its non-crossing diagonals  that are symmetric with respect to the axis. The poset    is isomorphic to the face poset of the associahedron $As_{k+1}$.
 \end{prop}
Proof. Cut the polygon through the symmetry axis  and contract the cut to a new vertex. This yields a  ($k+1$)-gon. Pairs of symmetric diagonals of the $2k$-gon correspond bijectively to diagonals of the ($k+1$)-gon. Single symmetric diagonals correspond to diagonals emanating from the new vertex.  We therefore arrive to a poset isomorphism. \qed

\medskip

The poset $\widetilde{D}^{inv}$ can be realized as the poset of some (uniquely defined) cell complex.
Indeed, let us build up  $\widetilde{\mathcal{D}}^{inv}$  starting from the cells of maximal dimension. Each such cell corresponds to cutting of the surface $F$ into a single polygon. Adding more diagonals reduces to Example \ref{ExSymmAssoc}. In other words,   $\widetilde{\mathcal{D}}^{inv}$ is a patch of associahedra.

As in \cite{GorPan}, and in Section \ref{SectMainConstr}, we  factorize $\widetilde{\mathcal{D}}^{inv}$ by the action of the pure mapping class group. For this purpose consider the defined below barycentric subdivision of $\widetilde{\mathcal{D}}^{inv}$.

\subsection*{Poset $\widetilde{BD}^{inv}$ and cell complex $\widetilde{\mathcal{BD}}^{inv}$.} We apply now the construction of the {order complex} of a poset, which gives us barycentric subdivision.
  Each  element of the  poset ${\widetilde{BD}}^{inv}_{g, b , n; n_1,... ,n_b}$ is (the strong equivalence class of) some symmetric admissible arrangement $A=\{d_1,...,d_m\} $ with a linearly ordered  partition $A=\bigsqcup S_i$   into some non-empty sets  $S_i$  such  that the first set $S_1$ in the partition  is an admissible arrangement, and all the sets $S_i$ are invariant under $inv$.


The partial order on $\widetilde{BD}^{inv}$ verbatim repeats the constructions of \cite{GorPan} and  Section \ref{SectMainConstr}. Namely, it is generated by the following rule:
 \newline $(S_1,...,S_p)\leq (S'_1,...,S'_{p'})$ whenever one of the two conditions holds:
\begin{enumerate}
  \item We have one and the same arrangement  $A$, and $(S'_1,...,S'_{p'})$ is an order preserving refinement of $(S_1,...,S_p)$.
  \item $p\leq p'$, and for all $i=1,2,...,p$, we have  $S_i=S'_i$. That is, $(S_1,...,S_p)$ is obtained from $(S'_1,...,S'_{p'})$ by removal
  $S'_{p+1},...,S'_{p'}$.
\end{enumerate}

By construction, the complex  $\widetilde{\mathcal{BD}}^{inv}$  is combinatorially isomorphic to the barycentric subdivision of $\widetilde{\mathcal{\mathcal{D}}}^{inv}$.

We are mainly interested in the  quotient complex:

\begin{definition}
  For a fixed data $(g, b , n; n_1,... ,n_b, \ inv)$, the diagonal  complex $\mathcal{BD}^{inv}_{g, b , n; n_1,... ,n_b}$ is defined as

  $$\mathcal{BD}^{inv}=\mathcal{BD}^{inv}_{g, b , n; n_1,... ,n_b}:=\mathcal{\widetilde{BD}}^{inv}_{g, b , n; n_1,... ,n_b}/PMC^{inv}(F).$$
  We define also
   $$\mathcal{D}^{inv}=\mathcal{{D}}^{inv}_{g, b , n; n_1,... ,n_b}:=\mathcal{\widetilde{D}}^{inv}_{g, b , n; n_1,... ,n_b}/PMC^{inv}(F).$$
\end{definition}
 Alternative definition reads as:

 \begin{definition}
   Each cell of the complex ${{\mathcal{BD}}}^{inv}_{g, b , n; n_1,... ,n_b}$ is labeled by the weak equivalence class of some symmetric admissible arrangement $A=\{d_1,...,d_m\} $ with a linearly ordered partition $A=\bigsqcup S_i$   into some non-empty sets  $S_i$  such  that the first set $S_1$   is an admissible arrangement, and all $S_i$ are invariant under $inv$.

The incidence rules  are the same as the above rules for the complex $\widetilde{\mathcal{BD}}^{inv}$.
 \end{definition}

\begin{prop} \label{2kgon} The cell complex $\mathcal{BD}^{inv}$ is regular. Its cells are combinatorial simplices.

\end{prop}
Proof.   If $(S_1,...,S_r)\leq (S'_1,...,S'_{r'})$ then there exists a unique (up to isotopy)  order-preserving  pure diffeomorphism of $F$ which  embeds \newline  $A= S_1\cup ...\cup S_r$ in  $A'=S'_1\cup ...\cup S'_{r'}$.
Indeed,   If $S_1=S_1'$,  the arrangement $S_1$  maps identically to itself since it has no automorphisms by stability assumption. The rest of the diagonals are diagonals in polygons, and are uniquely defined by their endpoints. Assume that $S_1 \subset S_1'$.
For the rest of the cases it suffices to take $A=S_1$, $A'=A =S_1'\bigsqcup S_2'$. If $A$ embeds in $A'$ in different ways,
then $A$ has a non-trivial isomorphism, which contradicts stability assumption.
\qed

\medskip

\medskip

\section{Symmetry vs hole}

In our setting, the surface $F$ is patched of two copies of an orientable surface with a distinguished boundary component $C^{fix}$. Contract $C^{fix}$ to a new vertex $v$ and
denote the resulted surface by $\frac{1}{2}F$. It inherits from the initial surface $F$ one half of its vertices.
Now replace the vertex $v$ by a hole. This gives $\frac{1}{2}F^{hole}$.

Take a simplex in ${\mathcal{BD}}^{inv}$  (starting from now, let us omit subscripts). It is labeled by  an admissible  arrangement $A=S_1\sqcup S_2...\sqcup S_k$.
Cut the surface $F$ with the arrangement through $C^{fix}$ and take one half. Contract $C^{fix}$ to a new vertex $v$. A moment's reflection reveals that one gets the surface $\frac{1}{2}F$ together with an admissible arrangement such that removal of $v$ with all the incident diagonals leaves an admissible arrangement.
Next, replace $v$ by a hole and remove  all the incident diagonals. We
obtain two well-defined maps:

  $${\mathcal{BD}}^{inv}\xrightarrow{\pi_1} \mathcal{BD}\Big(\frac{1}{2}F\Big),\hbox{and}$$
$$Im ~\pi_1 \xrightarrow{\pi_2} \mathcal{BD}\Big(\frac{1}{2}F^{hole}\Big).$$

\begin{thm}\label{ThmHalf} (1) The  map  $${\mathcal{BD}}^{inv}\xrightarrow {\pi_1} \mathcal{BD}\Big(\frac{1}{2}F\Big)$$ which takes one half of $F$ and contracts $C^{fix}$ to a new vertex $v$, is a combinatorial isomorphism on its image.

(2) The triple $$Im~ \pi_1 \xrightarrow{\pi_2} \mathcal{BD}\Big(\frac{1}{2}F^{hole}\Big)$$  is a homotopy equivalence. Here the map $\pi_2$ removes all diagonals incident to $v$ and replaces
$v$ by a hole.

(3) Altogether, $${\mathcal{BD}}^{inv}\xrightarrow{\pi_2 \circ \pi_1} \mathcal{BD}\Big(\frac{1}{2}F^{hole}\Big)$$  is a homotopy equivalence.
\end{thm}

Proof.

The claim (1) is clear by construction and Lemma \ref{Lemma1}.

 Assume that a simplex $\sigma$ belongs to $\mathcal{BD}\Big(\frac{1}{2}F\Big)$. As we know, it is labeled by some $(S_1,...,S_k)$, and the diagonals from $S_1$ cut $\frac{1}{2}F$ into disks.
The simplex $\sigma$ belongs to $Im ~\pi_1$ iff none of the discs has more than one corner incident to  $v$.

The statement  (3) follows directly from (1),(2).

\begin{figure}[h]
\centering \includegraphics[width=14cm]{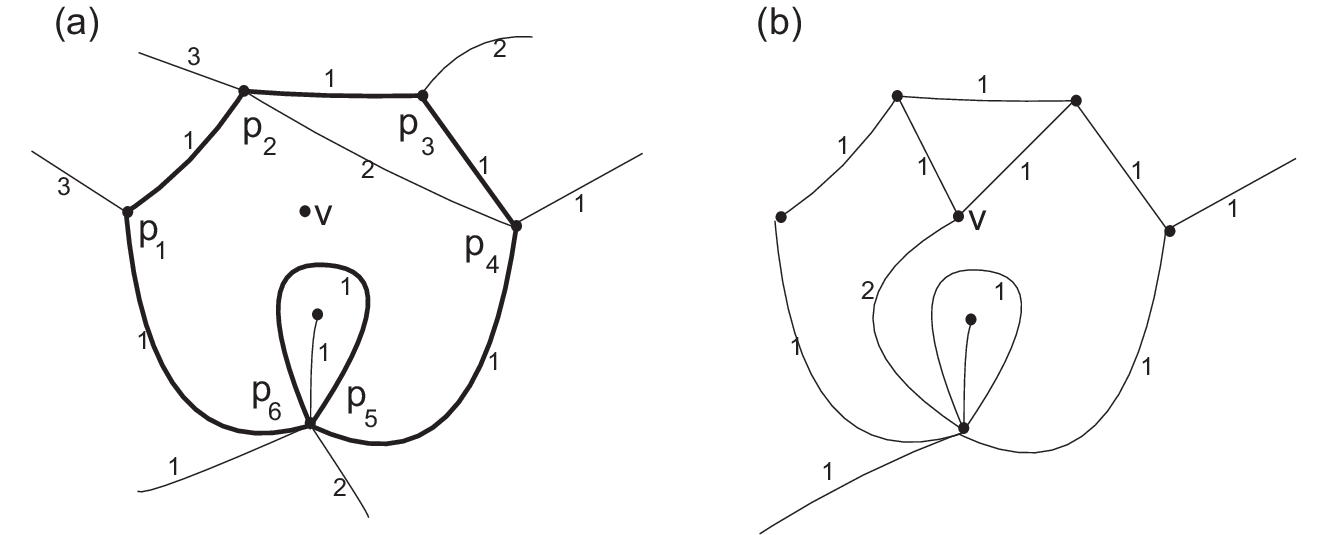}
\caption{(a) We depict a fragment of the label $(S_1,...,S_k)$ of one of the simplices $\sigma$. The polygon $P$ is the bold one. (b) The label of one of the simplices lying in the preimage $\pi_2^{-1}(Cl \sigma)$.}\label{PolygonP}
\end{figure}

Now prove (2).  Due to \cite{Q}, Theorem A,  it suffices to prove that the preimage of each closed simplex is contractible. So let us take a simplex $\sigma$ in $\mathcal{BD}\Big(\frac{1}{2}F^{hole}\Big)$. It is labeled by some $(S_1,...,S_k)$, where $S_1$ and   $A=S_1\sqcup...\sqcup S_k$ are admissible arangements.  Turn the hole to a vertex $v$. The vertex $v$ lies in the minimal polygon $P$ whose edges belong to $S_1$ {or are the edges of boundary components}. Denote the vertices  (with possible repetitions) of $P$  by  $p_1,...,p_q$, see Fig. \ref{PolygonP} (a).  A repeated vertex is given  different indices, although  it is one and the same marked point.

The labels of the simplices of the preimage
$\pi_1^{-1}(Cl\sigma)$ of the closure are obtained by: (1) eliminating some (might be none) of the entries from the end of the string  $(S_1,..., S_k)$, (2) replacing some consecutive  entries by their union, and (3) adding a number of \textit{new} diagonals in the polygon $P$ together with assigning numbers to the new diagonals. Each new diagonal should  emanate from $v$.  At least one of the new diagonals should be assigned number $1$, see Fig. \ref{PolygonP}(b).

\bigskip

We shall prove that $\pi_1^{-1}(Cl\sigma)$ is contractible by presenting a discrete Morse function with exactly one critical simplex (see Appendix A  for necessary backgrounds).

\bigskip

 \textbf{Step 1. "Moving new diagonals with numbers greater than $1$".}
Each new diagonal connects $v$ with some of  $p_k$. Since it is uniquely determined by $p_k$, let us call it $d_k$.
\begin{enumerate}
  \item  Assume that a simplex in the preimage $\pi_1^{-1}(Cl\sigma)$ is labeled by $(Q_1 ,..., (Q_i\cup \{d_1\}),..., Q_m)$  such that $i>1$ and $Q_i\neq \emptyset$.
Match it with $(Q_1 ,..., \{d_1\}, Q_i,..., Q_m)$.

The unmatched simplices are labeled by $(Q_1 ,..., Q_m )$ such that  either (a) $d_1\in Q_1$, or (b) $d_1$ is missing, or (c) $\{d_1\}$ is a singleton at the end.
  \item Take the simplices in $\pi_1^{-1}(Cl\sigma)$ that are unmatched on the previous step.
  Assume that a (nonmatched) simplex is labeled by $(Q_1 ,..., (Q_i\cup \{d_2\}),..., Q_m)$  such that $i>1$, $Q_i\neq \emptyset$, and $Q_i \neq \{d_1\}$.
Match it with $(Q_1 ,..., \{d_2\}, Q_i,..., Q_m)$.

  After this series of matchings the unmatched simplices have labels of the following four types: \begin{enumerate}
                                                                                  \item $(Q_1 ,..., Q_m )$ such that none of $Q_i$ contains $d_1$ or $d_2$ if $i>1$.
                                                                                  \item $(Q_1 ,..., Q_m,\{d_1\} )$ such that none of $Q_i$ contains  $d_2$ if $i>1$.
                                                                                      \item $(Q_1 ,..., Q_m,\{d_2\} )$ such that none of $Q_i$ contains  $d_1$ if $i>1$.
                                                                                  \item $(Q_1 ,..., Q_m,\{d_2\},\{d_1\} )$.
                                                                                \end{enumerate}

  \item Proceed the same way with all the other diagonals.
\end{enumerate}

Finally, the unmatched  simplices are labeled by $(Q_1 ,..., Q_m , \{d_{i_1}\}, ..., \{d_{i_r}\})$ such that (a)  none of $Q_i$ contains  new diagonals if $i>1$, and (b) $i_1>i_2>...>i_r$.

In other words,  for an unmatched simplex, the new diagonals either sit in the set $Q_1$, or come  as singletons at the very end of the label in the decreasing order.

\bigskip

 \textbf{Step 2. "Getting rid of old diagonals with numbers greater than $1$".}  Take a  simplex in  $\pi_1^{-1}(Cl\sigma)$  labeled by $(Q_1 ,..., Q_m , \{d_{i_1}\}, ..., \{d_{i_r}\})$ which is not matched on the Step 1. Denote the set of  all the new diagonals lying in $Q_1$  by $NEW$. Recall also that $Q_2,...,Q_m$ contain no new diagonals. Assuming that $Q_1\setminus NEW \neq S_1$, match  $(Q_1 ,..., Q_m , \{d_{i_1}\}, ..., \{d_{i_r}\})$ with  $(S_1\cup NEW, Q_1\setminus NEW \setminus S_1 , Q_2 , ..., Q_m , \{d_{i_1}\}, ..., \{d_{i_r}\})$.

 After Step 2, the unmatched simplices are labeled by $(Q_1 , \{d_{i_1}\}, ...,
  \{d_{i_r}\})$ such that (a)  $Q_1=S_1\cup NEW$, and (b)   $i_1>i_2>...>i_r$.

\bigskip

 \textbf{Step 3. "Final contractions. Getting rid of singletons with new diagonals."}
\begin{enumerate}
  \item Take an unmatched  (after Steps 1 and 2) simplex $(Q_1 , \{d_{i_1}\}, ...,
  \{d_{i_r}\})$. If $i_r=1$, match it with $(Q_1 , \{d_{i_1}\}, ...,
  \{d_{i_{r-1}}\})$. The unmatched simplices are those with $d_1\in Q_1$.
  \item Take an unmatched   simplex $(Q_1 , \{d_{i_1}\}, ...,
  \{d_{i_r}\}))$. If $i_r=2$, match it with $(Q_1 , \{d_{i_1}\}, ...,
  \{d_{i_{r-1}}\})$. The unmatched simplices are those without $d_1,d_2 \in Q_1$.
  \item Proceed the same way for $d_3,d_4,$ etc.
\end{enumerate}

Let us show that the matching is a discrete Morse function. The first two axioms follow straightforwardly from the construction.
The aciclicity will be proved a bit later.

So we arrive at  a discrete Morse function with a unique unmatched (that is, unique critical) simplex labeled by $(S_1\cup\{d_1,...,d_r\})$.
By the basic discrete Morse theory (see Appendix A), the preimage   $\pi_1^{-1}(Cl\sigma)$  is contractible.

\bigskip

Before we prove the acyclicity, let us look at an example of a gradient path:
$$(Q_1,Q_2\cup \{d_2\}, Q_3\cup \{d_1\}),\ \ (Q_1,Q_2\cup \{d_2\},\{d_1\}, Q_3),$$
$$  (Q_1,Q_2\cup \{d_2\}\cup\{d_1\}, Q_3), \ \ (Q_1,\{d_1\}, Q_2\cup \{d_2\}, Q_3)$$
$$  (Q_1\cup\{d_1\},Q_2\cup \{d_2\}, Q_3), \ \ (Q_1\cup\{d_1\}, \{d_2\}, Q_2, Q_3)$$

 Assume there exists a closed path.

(a) For a closed path, no new diagonal enters $Q_1$, since it can never leave $Q_1$.

(b) For each entry $\beta_i^{p+1}$ or $\alpha_{i}^p$  of the path, denote by $\overline{NEW}$  the set of new diagonals not lying in $Q_1$.
For a closed path, the set $\overline{NEW}$ cannot  decrease (if it contains some $d_i$, it never disappears on consequent steps of the path). Since the path is closed, $\overline{NEW}$ does not change during the path. Therefore, in a closed path the Step 3 is missing.

(c)
Assume that $d_i, d_j \in \overline{NEW}$  with $i<j$.   Then if $d_i$ is positioned to the left of $d_j$, it never appears to the right of $d_j$ in a closed path. We conclude that in a closed path,  all the entries of $\overline{NEW}$ appear as singletons coming in the decreasing order at the end of the label. In other words, in a closed path, Step 1  matchings are missing.

(d) Finally, a closed path cannot have Step 2 matchings only. This follows from a simple case analysis.
\qed

\section*{Appendix A. Discrete Morse theory \cite{Forman1}, \cite{Forman2}}

Assume we have a regular cell complex. By $\alpha^p, \ \beta^p$ we
denote its $p$-dimensional cells, or \textit{$p$-cells}, for short.

 A \textit{discrete vector field} is a set of pairs
$$\big(\alpha^p,\beta^{p+1}\big)$$
 such that:
\begin{enumerate}
    \item  each cell of the complex is matched with at most one
    other cell, and
    \item  in each pair, the cell $\alpha^p$ is a facet of $\beta^{p+1}$.

\end{enumerate}

Given a discrete vector field, a \textit{path}  is a sequence of
cells

$$\alpha_0^p, \ \beta_0^{p+1},\ \alpha_1^p,\ \beta_1^{p+1}, \ \alpha_2^p,\ \beta_2^{p+1} ,..., \alpha_m^p,\ \beta_m^{p+1},\ \alpha_{m+1}^p,$$
which satisfies the conditions:
\begin{enumerate}
    \item  Each $\big(\alpha_i^p$ and$\ \beta_i^{p+1}\big)$ are matched.
    \item Whenever $\alpha$ and  $\beta$ are neighbors in the path,
    $\alpha$ is a facet of $\beta$.
    \item $\alpha_i\neq \alpha_{i+1}$.
\end{enumerate}

A path is a \textit{closed path} if $\alpha_{m+1}^p=\alpha_{0}^p$.

A \textit{discrete Morse function on a regular cell complex} is a
discrete vector field without closed paths.
It gives a way of contracting all the cells of the complex that are matched:  if a cell $\sigma$ is matched with its facet\footnote{that is, a cell of dimension  $dim(\sigma)-1$ lying on the boundary of $\sigma$.}  $\sigma'$, then these two can be contracted by pushing $\sigma'$  inside $\sigma$.  Acyclicity  guarantees that if we have many matchings at a time, one can consequently perform the contractions. The order of contractions does not matter, and eventually one arrives at a complex homotopy equivalent to the initial one.

In the paper we use the following fact: \textit{if a regular cell complex has a discrete Morse function with exactly one critical cell, then
the complex is contractible.}

\section*{Appendix B. Associahedron and cyclohedron}

\subsection*{Associahedron and cyclohedron, \cite{Sta} and \cite{BottT}}  Assume that $n>2$ is fixed. We say that two diagonals in a convex $n$-gon are \textit{non-intersecting} if they  intersect  only at their endpoints (or do not intersect at all).  Consider all  collections of pairwise non-intersecting diagonals \footnote{It is important that the vertices of the polygon are labeled, and therefore
we do not identify collections of diagonals that differ on a rotation.} in the $n$-gon. This set is partially ordered by reverse inclusion,
and it was shown by John Milnor, that the poset is isomorphic to the face poset of some convex $(n-3)$-dimensional polytope $As_n$ called \textit{associahedron}.

In particular, the vertices of the associahedron $As_n$ correspond to the triangulations of the $n$-gon, and the edges correspond to edge flips in which one of the diagonals is removed  and replaced by a (uniquely defined) different diagonal.
Single diagonals  are in a bijection with facets of $As_n$, and  the empty set corresponds to the entire  $As_n$.

\medskip

 Fnalogously, centrally symmetric collections of diagonals in a $2n$-gon give rise to a convex polytope called \textit{cyclohedron},   or \textit{Bott–Taubes polytope}. Its first definition is the compactification of the configuration space of n points on the circle. 

\end{document}